\pdfoutput=1
\documentclass[manuscript,screen,nonacm]{acmart}
\AtBeginDocument{%
  }

\setcopyright{acmlicensed}
\copyrightyear{2026}
\acmYear{2026}
\acmDOI{XXXXXXX.XXXXXXX}

\acmJournal{TOMS}
\acmVolume{...}
\acmNumber{...}
\acmArticle{...}
\acmMonth{...}

\newcommand{\eg}{\emph{e.g.}}
\newcommand{\cl}{\operatorname{cl}}

\usepackage{minted}
\usepackage{tikz}
\usepackage{amsmath}
\usepackage{mathtools}
\theoremstyle{remark}
\newtheorem*{remark}{Remark}
\usepackage{subcaption}
\usepackage{placeins}

\begin{document}
\title{CVXPY 1.9: Recent Advances in Optimization Modeling Software}

\author{William Zhang}
\affiliation{%
  \institution{Polytechnique Montreal and Gridmatic}
  \city{Montreal}
  \country{Canada}
}
\email{william-2.zhang@polymtl.ca}

\author{Parth Nobel}
\affiliation{%
  \institution{Optimal Intellect}
  \city{San Francisco}
  \state{California}
  \country{USA}}
\email{ptnobel@alumni.stanford.edu}

\author{Aryaman Jeendgar}
\affiliation{%
  \institution{Technical University of Munich and International Computer Science Institute}
  \city{Heilbronn}
  \state{Baden-Württemberg}
  \country{Germany}
}
\email{aryaman.jeendgar@tum.de}

\author{Riley Murray}
\affiliation{%
  \institution{Sandia National Laboratories}
  \city{Livermore}
  \state{California}
  \country{USA}
}
\email{rjmurr@sandia.gov}

\author{Philipp Schiele}
\affiliation{%
  \institution{Citadel}
  \city{New York}
  \state{New York}
  \country{USA}
}
\email{phschiele@gmail.com}

\author{Steven Diamond}
\affiliation{%
  \institution{Optimal Intellect}
  \city{Mountain View}
  \state{California}
  \country{USA}
}
\email{diamond@cs.stanford.edu}

\renewcommand{\shortauthors}{Zhang et al.}
\newcommand\plotzero{{\color{gray}0}}
\begin{abstract}
CVXPY is a Python-embedded domain-specific language for convex optimization that lets users express problems in mathematical notation while the system verifies convexity and reduces valid programs to solver-ready form. This paper reports on the major advances from versions 1.1 through 1.9. These include a unified conic quadratic program (CQP) standard form for canonicalization; a stacked-slices backend that accelerates parameterized programs; first-class support for N-dimensional expressions; explicit sparsity for variables; support for multiple variable attributes; cones/atoms relevant to quantum information theory; and the introduction of disciplined nonlinear programming (DNLP). We outline the design, algorithms, and modeling consequences of these features.
\end{abstract}



\maketitle
\section{Introduction}
CVXPY began in 2013 as a Pythonic modular implementation of the modeling principles of CVX~\cite{grant2004dcp} in MATLAB.
A fundamental design decision of CVXPY is the separation between modeling, canonicalization, and solver interfaces~\cite{diamond2016cvxpy}.
Modeling in CVXPY is based on a grammar called disciplined convex programming (DCP), which allows for systematic verification of convexity.
CVXPY version 1.0 introduced a general reduction-based rewriting system for canonicalization~\cite{agrawal2018rewriting}, which enabled grammars such as disciplined geometric programming (DGP)~\cite{agrawal2019dgp} and disciplined quasiconvex programming (DQCP)~\cite{agrawal2020dqcp}, each being reduced to DCP-compatible forms.
CVXPY version 1.1 added disciplined parametrized programming (DPP)~\cite{agrawal2019differentiable}, which allows users to efficiently update canonicalized programs under a restricted grammar and enables differentiable programs and code generation (\texttt{CVXPYlayers}~\cite{agrawal2019differentiable}, \texttt{CVXPYgen}~\cite{cvxpygen,schaller2025codegenerationsolvingdifferentiating}).
Other grammars which focus on different problem classes have been implemented as separate standalone packages.
Examples include DCCP, DMCP, DCSP, and DSP~\cite{dccp, dmcp, dcsp, schiele2024dsp}.
CVXPY 1.9 added support for disciplined nonlinear programming (DNLP), a new syntax for specifying nonlinear programs~\cite{dnlp}.

The rest of the paper is structured as follows: \S 2 presents relevant standard forms while \S 3 presents canonicalization, the process of converting between these standard forms.
Finally, the paper concludes with a wide range of new canonicalization and modeling features in \S 4 and \S 5.

\section{Standard Forms}

CVXPY is a compiler between functional convex programming, a user-friendly representation of convex optimization problems, and conic quadratic programming (and its special cases), a solver-friendly representation of convex optimization problems.
Recognizing that users frequently want to solve similar problems, CVXPY supports parameterization, enabling reuse of compilation.

\subsection{Functional Convex Programming}

When a user specifies a problem in CVXPY, the first step is to verify that the problem is convex.
Users specify problems---by providing expression trees for each function---as the standard ``functional'' form for convex optimization
\begin{equation}\label{eq:functional-standard-form}
\begin{array}{ll}
\mbox{minimize} & f_0(x; \theta) \\
\mbox{subject to} & f_i(x; \theta) \leq 0, \ i=1,\ldots,m_1,\\
& g_i(x; \theta) = 0, \ i=1,\ldots,m_2,
\end{array}
\end{equation}
with variable $x\in\mathbf R^n$ and parameter vector $\theta\in\mathbf R^p$~\cite{Boyd_Vandenberghe_2004}.
When $\theta$ is fixed, we refer to \eqref{eq:functional-standard-form} as a problem instance, otherwise we refer to \eqref{eq:functional-standard-form} as a problem family.

Theorem~\ref{comp-theorem} provides a sufficient criterion for verifying convexity of the functional form.
The process of applying this convexity theorem on every node of the expression tree is known as disciplined convex programming (DCP).
It is disciplined in the sense that it restricts the choice of representation of mathematical expressions. 
For example $\sqrt{x^2}$ is not DCP, whereas $|x|$ is.
In addition to per-argument monotonicities, CVXPY uses sign information to extend the set of valid compositions.
Without signed DCP, $f(g(x)) = (|x|)^2$ is invalid since $f(t)=t^2$ is not monotone on $\mathbf{R}$.
With signed DCP, it is valid because $f$ is non-decreasing on $\mathbf{R}_{+}$ and the image of $g$ is $\mathbf{R}_{+}$.

\begin{theorem}[Convex Composition Rule]\label{comp-theorem}
\par $h(f_1(x), f_2(x), \dots, f_k(x))$ is convex when $h$ is convex and for each $i \in \{1, \dots, k\}$ the following holds:
\begin{itemize}
    \item $h$ is non-decreasing in argument $i$ and $f_i$ is convex, or
    \item $h$ is non-increasing in argument $i$ and $f_i$ is concave, or
    \item $f_i$ is affine.
\end{itemize}
The same results hold for concave functions by swapping any mention of convex with concave (and vice versa). Affine expressions require $h$ and all of its arguments to be affine~\cite{grant2004dcp}.
\end{theorem}

\subsection{Conic Quadratic Programming}

Quadratic programs generalize linear programs by allowing quadratic objective functions.
Cone programs generalize linear programs by abstracting constraints of the form
``$b - Ax \succeq 0$'' to ``$b - Ax \in \mathcal K$'', where $\mathcal K$ is a convex cone.\footnote{That is, $\mathcal K$ is a convex set such that $\theta x \in \mathcal K$ for all $x \in \mathcal K$ and all $\theta \geq 0$.}
We say \emph{conic quadratic program} in reference to a convex optimization problem with conic constraints and a linear or quadratic objective.

\begin{remark}
    Conic quadratic programming is not a major generalization of conic programming. If we allow the use of the second-order cone, then a conic quadratic program can have its objective linearized at the expense of one additional variable and one additional constraint with dimension $1$ plus the rank of the quadratic objective matrix.
    Nevertheless, retaining the quadratic objective can be useful from an algorithmic standpoint.
\end{remark}

Specialized solvers exist for conic quadratic programs~\cite{SCS3,CLARABEL,cuclarabel}, cone programs~\cite{domahidi2013ecos, mosek}, and quadratic programs~\cite{osqp}.
CVXPY interfaces with any such solver by first canonicalizing to a conic quadratic program
\begin{equation}
      \begin{array}{ll}
        \text{minimize} & \frac{1}{2} x^T P x + q^T x\\
        \text{subject to} &  b-Ax \in \mathcal K,\\
      \end{array}
      \label{eq:cp}
\end{equation}
where $x\in \mathbf{R}^n$ is the variable, $\mathcal{K}\subseteq \mathbf{R}^m$ is a closed convex cone, and the \emph{problem data} are $P \in \mathbf{S}^{n}_{+}$, $q \in \mathbf{R}^n$, $A \in \mathbf{R}^{m \times n}$, and $b\in \mathbf{R}^m$.
CVXPY supports problems where the cone is a product of the zero cone ($\{0\}$), the nonnegative cone ($\mathbf{R}_+$), the second-order cone ($\mathcal K_{\textrm{soc},n} = \{(t, x) \mid \|x\|_2 \leq t, t \in \mathbf{R}, x \in \mathbf{R}^{n-1}\}$), the power cone ($\mathcal K_{\textrm{pow},\alpha} = \{(x, y, z) \mid x\geq 0, y\geq 0, x^\alpha y^{1-\alpha} \geq |z|\}$), the exponential cone ($\mathcal K_{\textrm{exp}}=\cl\{(x, y, z) \mid y > 0, y\exp(x / y) \leq z\}$), and the cone of positive semidefinite matrices ($\mathbf{S}_+$).
Previously, CVXPY used two standard forms, cone programs, for which $P = 0$, and quadratic programs for which $\mathcal K$ is a product of only the zero cone and the nonnegative cone.
Further elaboration on unifying the two standard forms is in \S\ref{quad-conic}.

\subsection{Parameterized Conic Quadratic Programming}

In parameterized conic quadratic programming, $P, q, A, b$ in \eqref{eq:cp} are functions of a parameter $\theta$.
DPP is a grammar that guarantees parametrized convex programs can be reduced to affine-solver-affine (ASA) form.
ASA form is a representation of an optimization problem such that $P, q, A, b$ are affine functions of $\theta$ and the solution of the original problem is an affine function of the solution to the canonicalized problem.

DPP consists of the same composition rule as DCP but with additional restrictions on expressions involving parameters. 
A parameter is a symbolic constant with known properties (\emph{e.g.}, sign and sparsity) and unknown numeric value. 
The following canonicalization theorem shows that DPP-compliant programs can be reduced to ASA form.

\begin{theorem}[Canonicalization Map]\label{theo-canon}
\par The canonicalization map for a disciplined parametrized program can be represented by sparse tensors $Q \in \mathbf R^{n \times (n+1) \times (p+1)}$ and $R \in \mathbf{R}^{m \times (n+1) \times (p+1)}$, where $m$ is the dimension of the constraints.
Letting $\tilde{\theta} \in \mathbf{R}^{p+1}$ denote the concatenation of $\theta$ and the scalar $1$, the quadratic cone problem data can be obtained as 
\begin{equation}
    [P\quad q] = \sum_{i=1}^{p+1} Q[:,:,i]\cdot \tilde{\theta}_i,
\qquad\text{and}\qquad
    [A \quad b] = \sum_{i=1}^{p+1} R[:,:,i]\cdot \tilde{\theta}_i.
\end{equation}
\end{theorem}
The proof can be adapted from~\cite{agrawal2019differentiable}.

\section{Canonicalization}\label{canon}
Canonicalization converts a verified DCP/DPP model into a cone program via a sequence of reductions~\cite{grant2004dcp,agrawal2018rewriting}. The front end verifies convexity by parsing the user-specified expression tree.
Then, reductions---\emph{e.g.}, negating the objective of maximization problems and turning variable attributes into explicit constraints---are applied. The last symbolic reduction is to expand nonlinear atoms to their epigraph form, resulting in an affine expression tree.
Subsequently, a numerical reduction reduces the affine expression tree into the sparse tensors from Theorem~\ref{theo-canon}.
Finally, solver-specific transformations are applied to match the form required by the solvers' interfaces.
At each reduction, the inverse transformations are also recorded, allowing the solution of the standard form problem to be mapped back to the user-specified problem.
%

As an illustrative example, adapted from~\cite{agrawal2019differentiable},
consider the following optimization problem
\begin{equation}
    \begin{array}{ll}
        \text{minimize} & \|Fx-g\|_2 + \mu \|x\|_2\\
        \text{subject to} & x \succeq 0\\
    \end{array}
      \label{eq:example}
\end{equation}
with variable $x \in \mathbf{R}^n$ and parameters $F\in \mathbf{R}^{m\times n}$, $g\in \mathbf{R}^m$, and $\mu \geq 0$. 
This problem formulation given by the following CVXPY code conforms to the DPP grammar:
\begin{minted}{python}
x = cp.Variable(n)
F = cp.Parameter((m, n))
g = cp.Parameter(m)
mu = cp.Parameter(nonneg=True)
objective = cp.Minimize(cp.norm(F @ x - g) + mu * cp.norm(x))
prob = cp.Problem(objective, [x >= 0])
\end{minted}
For brevity, all code snippets omit the definition of dimension parameters (\eg, $n$ and $m$ in this example) and the imports of \texttt{cvxpy} and \texttt{numpy}.

In the following subsections we describe in greater detail the canonicalization procedure for problem~\eqref{eq:example}.
In the first phase, a front end takes a human-readable specification of a program and converts it to an intermediate representation;
both compilers and rewriting systems often use abstract syntax trees as their intermediate representations~\cite{dragonbook,agrawal2018rewriting}.

\subsection{Expression Tree Representation}
The first canonicalization step is to verify that the objective and constraint expressions of a CVXPY problem conform to the specified rule set.
If a problem does not conform to the rule set, CVXPY raises an error message.
This phase of the canonicalization is known as \textit{parsing}.
Throughout the parsing, an intermediate tree representation, also known as an \textit{abstract syntax tree} (AST), is formed. 
In the AST, each parent node represents an atomic function, called an \emph{atom}, whose children nodes are the arguments to the atom. 
The leaf nodes of the tree are variables, parameters, or constants. 
A DPP syntax tree for problem~\eqref{eq:example} is shown in Figure~\ref{dcp-syntax-tree}.
\newcommand\annotatedop[2]{\phantom{\ensuremath{#2}}\hspace{2ex}\ensuremath{#1}\hspace{2ex}{\color{gray}\ensuremath{#2}}}

\begin{figure}[ht]
\centering
\begin{subfigure}[t]{0.65\textwidth}
\centering
\begin{tikzpicture}
  [level 1/.style={sibling distance=40mm},
   level 2/.style={sibling distance=20mm}]
  \node { \annotatedop{+}{\|Fx-g\|_2 + \mu\|x\|_2}}
    child {node {\annotatedop{\|\cdot\|_2}{\|Fx-g\|_2}}
      child {node {\annotatedop{-}{Fx-g}}
          child {node {\annotatedop{\texttt{matmul}}{Fx}}
                child {node[draw] {$F$}}
                child {node[circle, draw] {$x$}} }
          child {node[draw] {$g$}}
            }
        }
    child {node {\annotatedop{\texttt{mul}}{\mu\|x\|_2}}
      child {node[draw] {$\mu$}}
      child {node {\annotatedop{\|\cdot\|_2}{\|x\|_2}}
            child {node[circle, draw] {$x$}}
      }
    };
\end{tikzpicture}
\caption{Abstract Syntax Tree of Objective}
\Description{A tree representing the DCP parsing of the objective}
\end{subfigure} %
\begin{subfigure}{0.3\textwidth}
\centering
\begin{tikzpicture}
  \node {\annotatedop{\succeq}{x \succeq 0}}
    child {node[circle, draw] {$x$}}
    child {node {$0$}};
\end{tikzpicture}
\caption{Abstract Syntax Tree of Constraints}
\Description{A tree representing the DCP parsing of the constraints}
\end{subfigure}
\caption{AST Representation. Users write expressions in Python code that is immediately parsed into an AST, shown here. DPP verification is then run on the AST to ensure canonicalization will be possible on the problem and to verify convexity. Variables are labeled with a circle, and parameters with a square. Leaf nodes that refer to the same parameter or variable are deduplicated in practice (not shown here).}
\label{dcp-syntax-tree}
\end{figure}
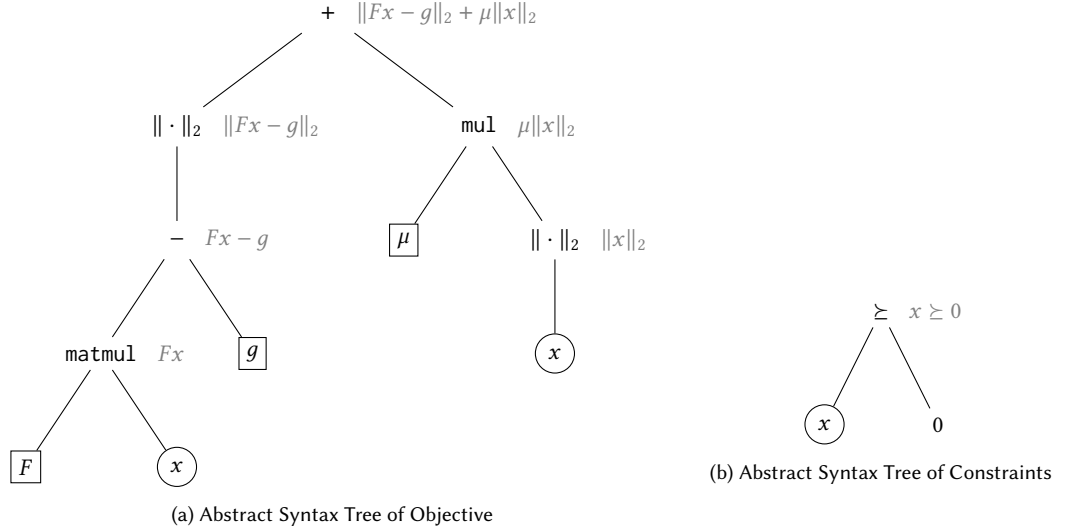

\subsection{Reductions}

Following the parsing and DPP verification, reductions are applied to transform the problem into a cone program.
In DPP, each atom has a \textit{canonicalization} that replaces nonlinear functions with an epigraph variable and auxiliary constraints of the form $f(x, t, \theta) \in \mathcal K$ where $t$ is the epigraph variable, $f$ is a biaffine function of the variables $(x, t)$ and the parameters $\theta$, and $\mathcal K$ is a convex cone.
One illustrative canonicalization replaces the $\|z_i\|_2$ expressions in the objective with a new epigraph variable $t_i$ and imposes the epigraph constraint in its conic form
\begin{equation}
(t_i, z_i) \in \mathcal K_{\textrm{soc}}.
\end{equation}

\subsection{Biaffine Operator Tree Representation}

\begin{figure}
\centering
\begin{subfigure}[t]{0.45\textwidth}
\centering
\begin{tikzpicture}
  [level 1/.style={sibling distance=20mm},
   level 2/.style={sibling distance=15mm}]
  \node {\annotatedop{+}{t_1+\mu  t_2}}
    child {node[circle,draw] {$t_1$}}
    child {node {\annotatedop{\texttt{mul}}{\mu  t_2}}
      child {node[draw] {$\mu$}}
      child {node[circle,draw] {$t_2$}}
    };
\end{tikzpicture}
\caption{Operator Tree of the Objective}
\Description{An operator tree representing the objective}
\end{subfigure} %
\begin{subfigure}[t]{0.5\textwidth}
\centering
\begin{tikzpicture}  
  [level 1/.style={sibling distance=20mm},
   level 2/.style={sibling distance=25mm},
   level 3/.style={sibling distance=15mm}]
  \node {\annotatedop{\texttt{vstack}}{(x,t_1,Fx-g,t_2,x)}}
    child {node[circle,draw] {$x$}}
    child {node[circle,draw] {$t_1$}}
    child {node {\annotatedop{+}{Fx - g}}
      child {node {\annotatedop{\texttt{matmul}}{Fx}}
        child {node[draw] {$F$}}
        child {node[circle,draw] {$x$}}
        }
      child {node {\annotatedop{\texttt{neg}}{-g}}
        child {node[draw] {$g$}}
        }
    }
    child {node[circle, draw] {$t_2$}}
    child {node[circle, draw] {$x$}};
\end{tikzpicture}
\caption{Biaffine Operator Tree of Constraints}
\Description{A linear operator tree representing the constraints}
\end{subfigure}
\caption{Biaffine Operator Trees. After the application of reductions, the constraints and the gradient with respect to the variables of the objective are biaffine functions. As with the AST representation, duplicate leaves are deduplicated in practice, which is not shown here. For readability, we show the objective rather than its gradient.}
\label{fig:lin-op-tree}
\end{figure}
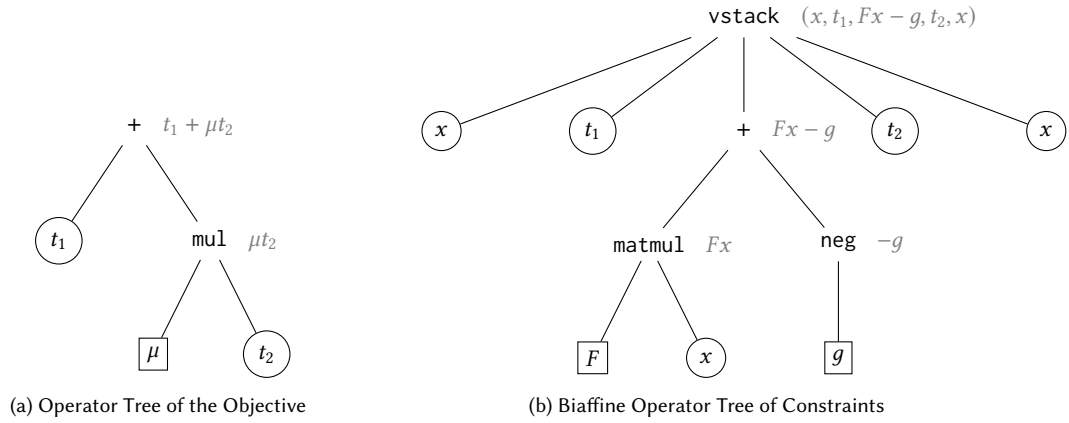

Following the reductions, the problem has been reduced to an equivalent parameterized convex cone program:
\begin{equation}\label{eq:graph}
\begin{array}{ll}
  \mbox{minimize}  & t_1 + \mu t_2 \\
  \mbox{subject to}  &
x \in \mathbf R_+^n \\ 
& (t_1, Fx -g) \in \mathcal K_{\mathrm{soc},m+1} \\ 
& (t_2, x) \in \mathcal K_{\mathrm{soc},n+1} \\ 
\end{array}
\end{equation}
with variables $t_1, t_2, x$.
We represent the biaffine function in the constraints and the gradient of the objective (which is an affine function since the objective is either quadratic or linear) as trees of biaffine operators, see Figure~\ref{fig:lin-op-tree}.

\subsection{Matrix Stuffing Representation}
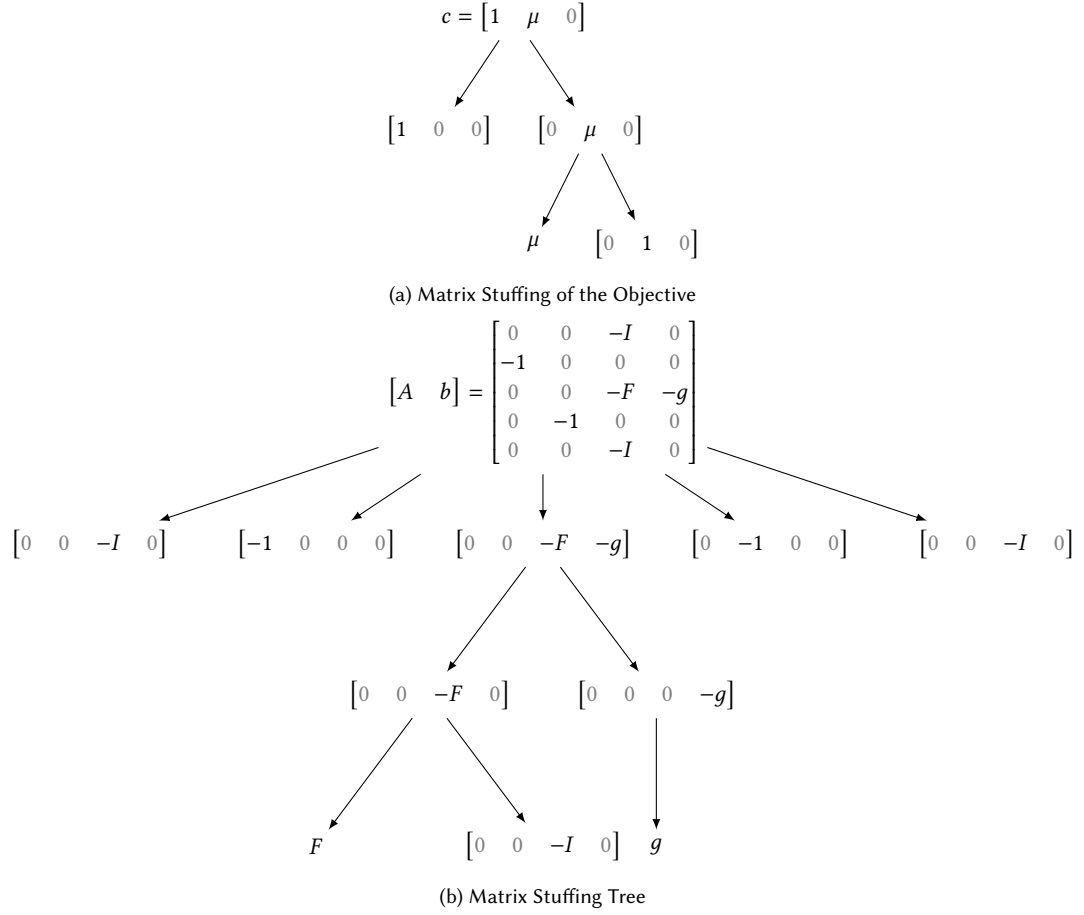
\begin{figure}
\begin{subfigure}{\textwidth}
\begin{center}
\begin{tikzpicture}
  [level 1/.style={sibling distance=20mm},
   level 2/.style={sibling distance=15mm},
    every node/.style={align=center},
    edge from parent/.style={draw, -latex}]
  \node {$c = \begin{bmatrix}
                  1 & 
                  \mu &
                  \plotzero
              \end{bmatrix}$}
    child {node {$\begin{bmatrix}
                  1 & 
                  \plotzero &
                  \plotzero
              \end{bmatrix}$}}
    child {node {$\begin{bmatrix}
                  \plotzero & 
                  \mu &
                  \plotzero
              \end{bmatrix}$}
      child {node {$\mu$}}
      child {node {$\begin{bmatrix}
                  \plotzero & 
                  1 &
                  \plotzero
              \end{bmatrix}$}}
      };
\end{tikzpicture}
\caption{Matrix Stuffing of the Objective}
\Description{Matrix stuffing of the objective}
\end{center}
\end{subfigure}
\begin{subfigure}{\textwidth}
\begin{center}
\begin{tikzpicture}[
    level distance=2cm,
    sibling distance=3cm,
    every node/.style={align=center},
    edge from parent/.style={draw, -latex}
]
\node {$\begin{bmatrix}
    A & b
\end{bmatrix} = \begin{bmatrix}
                   \plotzero & \plotzero & -I & \plotzero\\
                    -1 & \plotzero & \plotzero & \plotzero\\
                    \plotzero & \plotzero & -F & -g\\
                    \plotzero & -1 & \plotzero & \plotzero\\
                    \plotzero & \plotzero & -I & \plotzero\\
                \end{bmatrix}$}
    child {node {$\begin{bmatrix}
        \plotzero & \plotzero & -I & \plotzero\\
        \end{bmatrix}$}
    }
    child {node {$\begin{bmatrix}
        -1 & \plotzero & \plotzero & \plotzero\\
        \end{bmatrix}$}
    }
    child {node {$\begin{bmatrix}
                    \plotzero & \plotzero & -F & -g
                \end{bmatrix}$}
        child {node {$\begin{bmatrix}
                    \plotzero & \plotzero & -F & \plotzero
            \end{bmatrix}$}
            child {node {$F$}}
            child {node {$\begin{bmatrix}
                    \plotzero & \plotzero & -I & \plotzero\\
            \end{bmatrix}$}}
        }
        child {node {$\begin{bmatrix}
                    \plotzero & \plotzero & \plotzero & -g
            \end{bmatrix}$}
            child {node {$g$}}
            }
    }
    child {node {$\begin{bmatrix}
         \plotzero & -1 & \plotzero & \plotzero\\
        \end{bmatrix}$}
    }
    child {node {$\begin{bmatrix}
        \plotzero & \plotzero & -I & \plotzero\\
        \end{bmatrix}$}
    };
\end{tikzpicture}
\caption{Matrix Stuffing Tree}
\Description{A matrix stuffing tree}
\end{center}
\end{subfigure}
\caption{Matrix Stuffing Representation. Matrix stuffing generates a 3D tensor. To depict these 3D tensors here, we leave the parameter values symbolic. Compare this figure with Figure~\ref{fig:lin-op-tree}, which shows the operations being applied at each node.}
\label{fig:matrix-stuffing}
\end{figure}

Finally, the biaffine operator trees are materialized as sparse tensors.
This step is called \emph{matrix stuffing} and is done by the \emph{canonicalization backend}.
Canonicalization backends apply a depth-first walk on the tree to generate the sparse tensors $Q$ and $R$ from Theorem~\ref{theo-canon}.
See Figure~\ref{fig:matrix-stuffing}.

\section{New Canonicalization Features}

\subsection{Quadratic Cone Program Standard Form}\label{quad-conic}
When CVXPY was originally designed, it created two canonicalization pathways: one for cone programs,
\[
\begin{array}{ll}
\text{minimize} & c^T x \\
\text{subject to} & A x + s = b  \\
 & s \in \mathcal K,
\end{array}
\]
and one for quadratic programs,
\[
\begin{array}{ll}
\text{minimize} & \frac{1}{2} x^T P x + q^T x \\
\text{subject to} & \ell \preceq A x \preceq u .
\end{array}
\]
Recognizing the opportunity of recent solver advancements~\cite{SCS3,CLARABEL}, we combined our two standard forms into a single canonicalization pathway that generalizes both forms.

We now canonicalize to
\[
\begin{array}{ll}
\text{minimize} & \frac{1}{2} x^T P x + q^T x \\
\text{subject to} & Ax + s = b   \\
 & s \in \mathcal K.
\end{array}
\]

The reduction to a quadratic cone program is a modification of the reduction to a linear cone program. 
Instead of expanding all nonlinear atoms into an affine objective and conic constraints, quadratic-representable nonlinear atoms are expanded into a quadratic objective and conic constraints.

Most atoms with quadratic objectives in their expansion can be expressed purely as a quadratic term $x^TPx$, with no conic constraints.
A few functions have more interesting expansions, such as the huber function 
\[
\text{huber}(x, M) = \begin{cases}x^2 &|x| \leq M  \\2M|x| - M^2&|x| >M.\end{cases}
\]
The huber function can be rewritten as 
\[
\begin{array}{ll}
  \mbox{minimize}  & t_1^2 + 2M|t_2| \\
  \mbox{subject to}  & t_1 + t_2 = x,
\end{array}
\]
where $t_1$ and $t_2$ are auxiliary variables.

CVXPY users may still specify that a problem be canonicalized as a cone program with a linear objective rather than with a quadratic objective.

\subsection{Disciplined Nonlinear Programming}
In version 1.9, CVXPY allows users to specify nonlinear programs (NLP) of the form
\begin{equation}
      \begin{array}{ll}
        \text{minimize} & f(x)\\
        \text{subject to} &  c(x) = 0\\
        & \ell \leq x \leq u,
      \end{array}
      \label{eq:nlp-form}
\end{equation}
where $\ell \in \mathbf{R}^n$ and $u \in \mathbf{R}^n$ are given variable bounds and $f: \mathbf{R}^n \to \mathbf{R}$ and $c: \mathbf{R}^n \to \mathbf{R}^m$ are twice differentiable functions that can be nonconvex.
In addition, CVXPY supports combining nonsmooth convex and concave functions with smooth functions, subject to restrictions on how nonsmooth functions can be used. As in DCP, nonsmooth functions are relaxed in a lossless way to a form that can be handled by NLP solvers.
CVXPY currently supports the open-source solvers \texttt{IPOPT}~\cite{wachter2006} and \texttt{UNO}~\cite{Vanaret2025}, and the commercial solvers \texttt{KNITRO}~\cite{byrd2006knitro} and \texttt{COPT}~\cite{copt}.
For more details on the DNLP syntax, grammar rules, and canonicalization, see~\cite{dnlp}.

\subsection{Stacked-Slices Backend}\label{stacked-slices}
In the hopes of replacing our deprecated C++-based backend \texttt{CVXCORE}, we developed a \texttt{SciPy.sparse}-based backend.
To represent three-dimensional sparse tensors, the initial Python backend held per-parameter slices as lists of sparse matrices, which incurred Python-loop overhead for parameters of large size.
The new \emph{stacked-slices} architecture vertically concatenates slices and applies the linear maps via Kronecker identities.
For a tensor $Q \in \mathbf R^{a \times b \times (p+1)}$ whose slices correspond to the augmented parameter vector $\tilde{\theta}=(\theta,1)\in\mathbf R^{p+1}$, let $Q_i = Q[:,:,i] \in \mathbf R^{a \times b}$ and let $\mathcal S(Q) \in \mathbf R^{(p+1)a \times b}$ denote the vertical stack of the slices $Q_i$. If $L \in \mathbf R^{c \times a}$ is a linear map applied to each slice and $I_{p+1} \in \mathbf R^{(p+1)\times(p+1)}$ is the identity matrix over the parameter and constant slices, then
\[
\mathcal S(LQ) = (I_{p+1}\otimes L)\,\mathcal S(Q).
\]
This change removes most Python loops over lists of slices and speeds up canonicalization significantly in many cases. A continuous \texttt{asv} benchmark suite tracks performance on standard workloads in finance, SDP, and differentiable settings~\cite{asv}.

In version 1.8, a variant of the stacked-slices backend was introduced that uses a custom N-dimensional sparse COO representation in place of \texttt{SciPy.sparse}.
The COO representation is more memory efficient, particularly for problems with large numbers of parameters.
CVXPY defaults to the COO backend for DPP problems with more than 1000 parameters.

\section{New Modeling Features}
This section covers recent improvements to the modeling abilities of CVXPY: support for N-dimensional expressions, a sparsity attribute, and multiple attributes.
This section also reviews CVXPY's features for modeling in quantum information theory; the implementation of these features also enabled approximate solution of mixed-integer exponential-cone programs (and mixed-integer geometric programs) with mixed-integer SOCP solvers.

\subsection{N-Dimensional Expressions}
Starting from version 1.6, CVXPY began supporting N-dimensional (N-d) expressions, intended to match the syntax of \texttt{ndarray}s from NumPy~\cite{Harris_2020}. 
N-dimensional data structures---also referred to as \textit{tensors}---are used widely across the scientific Python ecosystem (\eg, \texttt{xarray}~\cite{Hoyer-2017}, \texttt{PyTorch}~\cite{paszke2019pytorchimperativestylehighperformance}, and \texttt{JAX}~\cite{jax2018github}).
Previously, expressions involving more than two dimensions had to be implemented by hand, either through nested lists of variables or custom indexing logic for a stacked variable.

As the standard form requires a single vector variable, N-dimensional expressions must be flattened by the backend. CVXPY defaults internally to the Fortran order, or column-major order, for flattening higher-dimensional arrays.
This process is identical to how we previously handled matrices.


\subsubsection{Sum Operator}

\paragraph{Example}
Consider an expression $\mathcal X$ with shape $(2,2,2)$ being summed over a single axis $s=\{1\}$.
The user considers $\mathcal X$ as a vector of matrices:
\[
\mathcal{X} = \begin{bmatrix}
\begin{bmatrix}
x_{111} & x_{112} \\
x_{121} & x_{122}
\end{bmatrix},
\begin{bmatrix}
x_{211} & x_{212} \\
x_{221} & x_{222}
\end{bmatrix}
\end{bmatrix}
\]
The user expects to see the output of the sum operator $\mathcal Y$ to be a matrix:
\[
\mathcal{Y} =
\begin{bmatrix}
y_{11} & y_{12} \\
y_{21} & y_{22} \\
\end{bmatrix}
=
\begin{bmatrix}
x_{111} + x_{121} & x_{112} + x_{122} \\
x_{211} + x_{221} & x_{212} + x_{222}
\end{bmatrix}.
\]
In the CVXPY implementation, both expressions are stored in their vectorized form and the sum operator reduces to the following sparse matrix:
\[
\begin{array}{c|cccccccc}
& x_{111} & x_{211} & x_{121} & x_{221} & x_{112} & x_{212} & x_{122} & x_{222} \\
\hline
y_{11} & 1 & \plotzero & 1 & \plotzero & \plotzero & \plotzero & \plotzero & \plotzero \\
y_{21} & \plotzero & 1 & \plotzero & 1 & \plotzero & \plotzero & \plotzero & \plotzero \\
y_{12} & \plotzero & \plotzero & \plotzero & \plotzero & 1 & \plotzero & 1 & \plotzero \\
y_{22} & \plotzero & \plotzero & \plotzero & \plotzero & \plotzero & 1 & \plotzero & 1 
\end{array}
\]



\subsubsection{Indexing Operator}

\paragraph{Example}
Consider an arbitrary expression with shape $(3,2,2)$ and the slice operation $\mathcal{X}[0:2, :, 1]$, which selects the first two slices along axis 0, all elements along axis 1, and index 1 along axis 2. First, we present the internal representation of the expression.
The original expression is of the following form:
\[
\mathcal{X} = \begin{bmatrix}
\begin{bmatrix}
x_{111} & x_{112} \\
x_{121} & x_{122}
\end{bmatrix},
\begin{bmatrix}
x_{211} & x_{212} \\
x_{221} & x_{222}
\end{bmatrix},
\begin{bmatrix}
x_{311} & x_{312} \\
x_{321} & x_{322}
\end{bmatrix}
\end{bmatrix}.
\]

The indexing operation $\mathcal{X}[0{:}2,~:~,~1]$ produces a tensor $\mathcal{Y}$ of shape $(2,2)$:
\begin{align*}
\mathcal{Y} =
\begin{bmatrix}
y_{11} & y_{12} \\
y_{21} & y_{22}
\end{bmatrix}
=
\begin{bmatrix}
x_{112} & x_{122} \\
x_{212} & x_{222}
\end{bmatrix}.
\end{align*}

Internally, CVXPY represents the indexing operation as a sparse matrix that maps the vectorization of the input tensor to the vectorization of the selected subtensor:
\begin{align*}
\begin{array}{c|cccccccc}
& x_{111} & \cdots & x_{112} & x_{212} & x_{312} & x_{122} & x_{222} & x_{322} \\
\hline
y_{11} & \plotzero & \plotzero & 1 & \plotzero & \plotzero & \plotzero & \plotzero & \plotzero \\
y_{21} & \plotzero & \plotzero & \plotzero & 1 & \plotzero & \plotzero & \plotzero & \plotzero \\
y_{12} & \plotzero & \plotzero & \plotzero & \plotzero & \plotzero & 1 & \plotzero & \plotzero \\
y_{22} & \plotzero & \plotzero & \plotzero & \plotzero & \plotzero & \plotzero & 1 & \plotzero
\end{array}.
\end{align*}

%

\subsection{Sparsity Attribute} 
In version 1.6, CVXPY introduced a sparsity attribute for variables, which allows users to specify the nonzero entries of a variable.
The sparsity attribute avoids storing unnecessary variables and provides performance improvements in both memory use and computation time during the canonicalization and solve steps.
The syntax for specifying the sparse indices is the same as SciPy's \texttt{sparse.coo\_array} format~\cite{Virtanen_2020}. 

\paragraph{Example}
Consider an arbitrary expression with shape $(2,2,3)$ that has only three nonzero entries, at indices $(0,1,0)$, $(1,0,1)$, and $(1,1,2)$. The user sees an expression of the form
\begin{align*}
\begin{bmatrix}
\begin{bmatrix}
\plotzero & \plotzero & \plotzero \\
z_1 & \plotzero & \plotzero
\end{bmatrix},
\begin{bmatrix}
\plotzero & z_2 & \plotzero \\
\plotzero & \plotzero & z_3
\end{bmatrix}
\end{bmatrix}
\end{align*}

Using the sparsity attribute, CVXPY can internally represent this expression with a compact vector of three entries $\mathcal{Z} = (z_1, z_2, z_3)$.
When instantiating the leaves of the matrix stuffing tree, instead of using an identity matrix to represent the $\mathcal Z$ tensor, we use the following sparse matrix:
\[
\begin{array}{c|ccc}
& z_1 & z_2 & z_3 \\
\hline
x_{000} & \plotzero & \plotzero & \plotzero \\
x_{100} & \plotzero & \plotzero & \plotzero \\
x_{010} & 1 & \plotzero & \plotzero \\
x_{110} & \plotzero & \plotzero & \plotzero \\
x_{001} & \plotzero & \plotzero & \plotzero \\
x_{101} & \plotzero & 1 & \plotzero \\
\vdots & \plotzero & \plotzero & \plotzero \\
x_{112} & \plotzero & \plotzero & 1
\end{array}
\]

%

\subsubsection{Defining a Sparse Variable}
\begin{minted}{python}
# Create an upper triangular sparse variable
X = cp.Variable((10, 10), sparsity=np.triu_indices(n=10))
\end{minted}

Another way to define a sparse variable is to use \texttt{np.where} on problem-specific data.

\begin{minted}{python}
# Define problem data (adapt to your use case)
data = np.random.randn(10, 10)
# Create a sparse variable given condition on data
x = cp.Variable((10, 10), sparsity=np.where(data > 0.5))
\end{minted}

Finally, you can also define the sparsity attribute manually. 
The input to the attribute needs to conform to the index format of \texttt{numpy.where} and \texttt{numpy.nonzero}.

\begin{minted}{python}
# Create a sparse variable manually
# The two tuples correspond respectively to row and column indices
x = cp.Variable((3, 3), sparsity=[(0, 1), (0, 2)])
# This is equivalent to calling np.where(data == 1) on the following matrix
# data = np.array([[1, 0, 0],
#                  [0, 0, 1],
#                  [0, 0, 0]])
\end{minted}

\subsection{Multiple Attributes}
In version 1.7, CVXPY began supporting multiple attributes for variables. 
CVXPY's variable attributes can be viewed as implicit constraints that provide useful sign and monotonicity information to the DCP analyzer. 
In addition, due to the large potential computational cost, projections for validating parameter inputs with multiple attributes are disabled. 
The following examples show some combinations of variable attributes that are supported in CVXPY 1.9. 

\subsubsection{Sparse Integer Variable}
\begin{minted}{python}
x = cp.Variable((2, 2), sparsity=[(0, 1), (0, 1)], integer=True)
prob = cp.Problem(cp.Minimize(cp.sum(x)), [x >= -5.5])
prob.solve()
np.testing.assert_allclose(prob.value, -10)
# Although the lower bound on x is -5.5, the nonzero entries
# of x are set to -5 due to the integer attribute
# .value is a dense NumPy array; .value_sparse returns a COO array
np.testing.assert_allclose(x.value, np.array([[-5, 0], [0, -5]]))
\end{minted}
\subsubsection{Bounded Integer Variable}
\begin{minted}{python}
x = cp.Variable((2, 2), integer=True, bounds=[-1.5, 2])
prob = cp.Problem(cp.Minimize(cp.sum(x)), [])
prob.solve()
np.testing.assert_allclose(prob.value, -4)
# Although the lower bound is -1.5, each entry of x
# is set to -1 due to the integer attribute
np.testing.assert_allclose(x.value, np.array([[-1, -1], [-1, -1]]))
\end{minted}
\subsubsection{Bounded Sparse Variable}
\begin{minted}{python}
x = cp.Variable((3, 3), sparsity=[(0, 1), (1, 2)], bounds=[-3.5, 7])
prob = cp.Problem(cp.Minimize(cp.sum(x)), [])
prob.solve()
np.testing.assert_allclose(prob.value, -7)
# Each nonzero entry of x is set to the lower bound of -3.5
np.testing.assert_allclose(x.value, np.array([[0, -3.5, 0],
                                              [0, 0, -3.5],
                                              [0, 0, 0.0]]))
\end{minted}

\subsection{Quantum Information Modeling}\label{sec:qi-modeling}

CVXPY has added support for some key cones and atoms relevant for modeling problems arising in quantum information theory, such as calculating channel capacities or entanglement measures. These can broadly be categorized in two groups.
 
\paragraph{Linear Algebra and Tensor Product Spaces.}

Tensor product spaces of complex matrices are ubiquitous in quantum information theory.
The familiar Kronecker product (available in CVXPY as the \texttt{kron} atom) is used to construct a matrix in a tensor product space $W = U \otimes V$ given matrices in the generating spaces $U$ and $V$. 

When $U$ and $V$ are spaces of square matrices we arrive at important generalizations of the matrix trace and the matrix transpose.
These generalizations can be defined in terms of orthonormal bases $B_U$ and $B_V$ for $U$ and $V$, which generate the basis $B_W = \{ u \otimes v \,:\; u \in B_U, v \in B_V\}$ for $W$.
\begin{itemize}
    \item The \emph{partial trace} with respect to $U$ is the function $\operatorname{Tr}_U : W \to V$ that maps $u \otimes v \in B_W$ to $\operatorname{Tr}(u) v$, and is then extended to all of $W$ by $\mathbf{C}$-linearity; the partial trace with respect to $V$ is the operator $\operatorname{Tr}_V : W \to U$ where $\operatorname{Tr}_V(u \otimes v) = \operatorname{Tr}(v) u$, again extended to all of $W$ by linearity.
    \item The \emph{partial transpose} (with respect to $U$ or $V$) is defined analogously.
    We emphasize that the partial transpose must use the \emph{transpose}, not the \emph{conjugate transpose}.
    The conjugate transpose is not a $\mathbf{C}$-linear operation.
\end{itemize}

CVXPY does not provide a first-class abstraction for vector spaces (tensor product, or otherwise).
One implicitly works with tensor product spaces simply by keeping a list of dimensions of the generating spaces.
In the context of quantum information, the generating spaces are called \emph{subsystems}.

\paragraph{Quantum Information Atoms.}
These atoms implement key concave and convex functions from quantum information theory, which often form the objective function of an optimization problem.
Matrices provided to these atoms must be Hermitian positive semidefinite (PSD).
\begin{enumerate}
    \item \texttt{von\_neumann\_entr}. The von Neumann entropy, $H(\rho) = -\operatorname{Tr}\left(\rho\log\rho\right)$, is the quantum mechanical analogue of Shannon entropy. The \texttt{von\_neumann\_entr} atom is concave.
    \item \texttt{quantum\_rel\_entr}. The quantum relative entropy, $D(\rho || \tau)=\operatorname{Tr}\left(\rho\left(\log\rho-\log\tau\right)\right)$, is the quantum analogue of the Kullback-Leibler (KL) divergence. It is a non-negative measure of \textit{distinguishability} between two quantum states $\rho, \tau$ and $D(\rho||\tau)=0$ only if $\rho=\tau$. The quantum relative entropy is jointly convex in its arguments $(\rho, \tau)$.
    \item \texttt{quantum\_cond\_entr}. Quantum conditional entropy is defined for a state $\rho \in U \otimes V$ and a subsystem $U$ or $V$. Conditioning with respect to $V$, it is equal to $H(\rho) - H(\operatorname{Tr}_U(\rho))$. This is concave in $\rho$, but is not naturally DCP. Representing its hypograph requires quantum relative entropy.
\end{enumerate}

We now demonstrate two quantum information problems that can be modeled with combinations of the atoms described above. These problems are borrowed from~\cite{FF18}.

\paragraph{\textbf{Relative Entropy of Entanglement} (REE)}
The state of a two-qubit quantum register is described by a trace-1 Hermitian PSD matrix (a \emph{density matrix}) of order $4$.
Such a state $\rho$ is called \emph{separable} if it can be expressed as a convex combination of tensor product states, $\rho_0 \otimes \rho_1$, where $(\rho_0, \rho_1)$ are density matrices of order $2$.
States that are not separable are called \emph{entangled}.
The REE of a state $\rho$ is the minimum quantum relative entropy achievable between $\rho$ and a separable state $\tau$ when we optimize over $\tau$.
For a two-qubit register, we can compute $\rho$'s REE by optimizing over $\tau$ that satisfy the \emph{positive partial transpose} condition.
It makes no difference whether the partial transpose is over the first or second qubit's subsystem.

\begin{listing}[htbp]
\begin{minted}{python}
# Define hypothetical input to the problem
U = np.random.randn(4, 4) + 1j * np.random.randn(4, 4)
psd = U @ U.T.conj()
rho = psd / np.trace(psd)
# Build and solve the problem
tau = cp.Variable((4, 4), hermitian=True)
func = cp.quantum_rel_entr(rho, tau) / np.log(2)
cons = [cp.partial_transpose(tau, [2, 2], 1) >> 0]
cons += [tau >> 0, cp.trace(tau) == 1]
prob = cp.Problem(cp.Minimize(func), cons)
prob.solve()
\end{minted}
\caption{Computing the Relative Entropy of Entanglement (REE) in CVXPY}
\label{lst:ree}
\end{listing}

\FloatBarrier

\paragraph{\textbf{Entanglement-Assisted Classical Capacity}.}
We try to calculate the fastest rate of transferring classical information (bits, 0s \& 1s) over a noisy quantum channel $\Phi$ when both the sender and receiver share an entangled state.
This is quantified by using the \textit{mutual information}, $I(\rho, \Phi)$, which can be implemented within CVXPY by a combination of the atoms \texttt{hermitian\_wrap}, \texttt{quantum\_cond\_entr}, \texttt{von\_neumann\_entr}, and \texttt{partial\_trace}.

\begin{listing}[htbp]
\begin{minted}{python}
# Define Phi by Stinespring rep., Phi(rho) = Tr_2(C rho C^\dagger)
C = np.array([
    [1, 0],
    [0, np.sqrt(0.01)],
    [0, np.sqrt(0.99)],
    [0, 0],
])
# Build and solve the CVXPY model
rho       = cp.Variable((2, 2), hermitian=True)
sigma     = cp.hermitian_wrap( C @ rho @ C.T.conj() )
Phi_rho   = cp.hermitian_wrap(cp.partial_trace(sigma, [2, 2], 1))
I_rho_Phi = cp.quantum_cond_entr(sigma, [2, 2], 0)
I_rho_Phi = (I_rho_Phi + cp.von_neumann_entr(Phi_rho)) / np.log(2)
obj = cp.Maximize(I_rho_Phi)
prob = cp.Problem(obj, [rho >> 0, cp.trace(rho) == 1])
prob.solve()
\end{minted}
\caption{Entanglement-Assisted Classical Capacity of a Quantum Channel in CVXPY}
\label{lst:en-classical-quantum-cap}
\end{listing}

\FloatBarrier

\subsection{Approximating Logs and Exponentials}
In \S\ref{sec:qi-modeling}, we outlined the various tools that CVXPY offers for modeling problems in quantum information theory.
This subsection discusses the machinery behind these atoms and points out some compelling corollary uses.

\paragraph{Operator Relative Entropy} The operator relative entropy cone for order-$n$ Hermitian matrices, $\mathbf{H}^n$, is given by
\[
K^n_{\text{re}} = \cl\left\{ (X,Y,T) \in \mathbf{H}^n_{++}\times\mathbf{H}^n_{++}\times\mathbf{H}^n \;\middle|\; \sqrt{X}\log(\sqrt{X} Y^{-1} \sqrt{X}) \sqrt{X} \preceq T  \right\}.
\]
The square roots and log in the display above are understood in the sense of spectral functions.

CVXPY's \texttt{OpRelEntrConeQuad} class provides approximate representations of $K^n_{\text{re}}$ using the techniques in the foundational work by Fawzi, Saunderson, and Parrilo~\cite{FSP18}.
The approximations, denoted $K^n_{m,k}$ for integer parameters $m,k$, are constructed using operator monotonicity of certain rational functions and Gauss-Legendre quadrature.
The main representation result,~\cite[Theorem 3]{FSP18}, requires $m + k + 1$ auxiliary Hermitian matrices and $(m + k)$ PSD constraints of sizes $2n\times 2n$ each.
The parameter $m$ controls the number of quadrature nodes and $k$ controls the number of iterative square roots applied to the argument to improve the approximation's accuracy, especially for values far from $1$.
The use of these approximate cones for quantum information modeling is described in~\cite{FF18}.

\paragraph{\texttt{RelEntrConeQuad}.}
The operator relative entropy approximation scheme when $n=1$ provides a way to approximate the exponential cone.
We introduced \texttt{RelEntrConeQuad} to handle this special case, with utility functions to map to and from the native \texttt{ExpCone}.
The $2 \times 2$ PSD constraints for \texttt{RelEntrConeQuad} are canonicalized into equivalent second-order cone constraints. 

This canonicalization has two potential uses.
The first is to open up the possibility of using SOCP solvers in place of less common exponential cone solvers.
This is particularly valuable for mixed-integer problems, in light of open-source mixed-integer SOCP solvers such as SCIP~\cite{SCIP}.
The second potential use is in improving convergence for solvers that support the exponential cone directly.

\paragraph{Market Equilibrium with Only Second-Order Cones}
In Listing~\ref{lst:market-equilibrium}, we show a market equilibrium problem with continuously divisible goods and linear utilities.
This problem can be formulated as an exponential cone program~\cite{Eisenberg1961}.
The code below builds the conic form directly, and shows how one can elect to 
replace the exponential cone with the \texttt{RelEntrConeQuad} approximation.
The approximate form of the problem is more tractable for many open-source solvers, including Clarabel.

\begin{listing}
\begin{minted}{python}
# Valuations (V), budgets (b), and utilities.
buyers, items = 400, 200
V = 10 ** np.random.uniform(-2, 2, (buyers, items))
b = 10 ** np.random.uniform(-2, 1, buyers)
b /= np.sum(b)
allocations = cp.Variable((buyers, items), nonneg=True)
utilities = cp.sum(cp.multiply(V, allocations), axis=1)
# Maximize conic form of Nash Social Welfare
aux_var = cp.Variable(buyers)
aux_con = cp.ExpCone(aux_var, np.ones(buyers), utilities)
objective = cp.Maximize(b @ aux_var)
constraints = [
    aux_con.as_quad_approx(m=3, k=3),  # Try dropping "as_quad_approx"
    cp.sum(allocations, axis=0) <= 1,
]
prob = cp.Problem(objective, constraints)
prob.solve(solver='CLARABEL', verbose=True)
\end{minted}
\caption{Approximating an Exponential Cone Program to Help Solve a Market Equilibrium Problem}
\label{lst:market-equilibrium}
\end{listing}

\section{Conclusion}

CVXPY has sought to improve its performance, to have its many features be increasingly cross-compatible, and to expand its modeling capabilities for new domains.
We aim to finish the deprecation and removal of \texttt{CVXCORE} as we seek to further improve performance.
We hope to develop a stable serialization format and to improve our infeasibility debugging support.
Our sister library CVXPYlayers has greatly improved and expanded our GPU support in ways that will further the applicability of CVXPY across domains and fields.

\FloatBarrier

\begin{acks}
We thank Gridmatic for funding the development of N-dimensional arrays.
We thank Google for Google Summer of Code.
We thank NumFOCUS for facilitating our applications to GSoC as well as a developer fund to improve CVXPY's website. 
Finally, we would like to thank the numerous CVXPY contributors and solver interface maintainers. 

Parth Nobel was supported in part by the National Science Foundation Graduate Research Fellowship
Program under Grant No. DGE-1656518. Any opinions, findings, and conclusions or recommendations
expressed in this material are those of the author(s) and do not necessarily reflect
the views of the National Science Foundation.
\end{acks}

\bibliographystyle{ACM-Reference-Format}
\bibliography{cvxpy-paper}

\end{document}